\documentclass[11pt, twoside]{amsart}
\usepackage{amsmath, amsthm, amscd, amsfonts, amssymb, graphicx, color}
\usepackage[bookmarksnumbered, plainpages, backref]{hyperref}

\begin{document}

\title{The Method of Almost Convergence with Operator of the form Fractional Order and Applications}

\author{Murat Kiri\c{s}ci*, U\v{g}ur Kadak}

\address{[Murat Kiri\c{s}ci]Department of Mathematical Education, Hasan Ali Y\"{u}cel Education Faculty,
Istanbul University, Vefa, 34470, Fatih, Istanbul, Turkey \vskip 0.1cm }
\email{mkirisci@hotmail.com, murat.kirisci@istanbul.edu.tr}

\address{[U\v{g}ur Kadak]Department of Mathematics,
Bozok University,Yozgat, Turkey \vskip 0.1cm }
\email{ugurkadak@gmail.com}

\thanks{*Corresponding author.}

\begin{abstract}
The purpose of this paper is twofold. First, basic concepts such as Gamma function, almost convergence,
fractional order difference operator and sequence spaces are given as a survey character. Thus, the current knowledge
about those concepts are presented. Second, we construct
the almost convergent spaces with fractional order difference operator and compute dual spaces which are help us in the characterization of
matrix mappings. After we characterize to the matrix transformations, we give some examples. In this paper, the notation $\Gamma(n)$ will be
shown the Gamma function. For $n\not \in \{0,-1,-2,\ldots\}$, Gamma function defined by an improper integral $\Gamma(n)=\int_{0}^{\infty}e^{-t}t_{n-1}dt$.
\end{abstract}

\maketitle

\pagestyle{plain} \makeatletter
\theoremstyle{plain}
\newtheorem{thm}{Theorem}[section]
\numberwithin{equation}{section}
\numberwithin{figure}{section}  
\theoremstyle{plain}
\newtheorem{pr}[thm]{Proposition}
\theoremstyle{plain}
\newtheorem{exmp}[thm]{Example}
\theoremstyle{plain}
\newtheorem{cor}[thm]{Corollary} 
\theoremstyle{plain}
\newtheorem{defin}[thm]{Definition}
\theoremstyle{plain}
\newtheorem{lem}[thm]{Lemma} 
\theoremstyle{plain}
\newtheorem{rem}[thm]{Remark}
\numberwithin{equation}{section}

\section{Introduction}

This paper consist of 5 Section. In section 1, we give some basic definitions of the Gamma functions,
almost convergent sequence spaces, fractional order difference operators and recall some basic knowledge
related to sequence spaces, Schauder basis, $\alpha-$, $\beta-$ and $\gamma-$ duals, matrix transformations.
Section 2 deal with the matrix domain of almost convergent space. Also, we compute dual spaces of new almost
convergent space, in Section 2. Section 3 is devoted to the matrix transformations between new space and classical spaces.
In section 4, we give some examples.
In last section, subsequent to giving a short analysis on the basic results of the present paper, some further suggestions are noted.\\

\subsection{Gamma Function}
In view of the analogy between infinite series and improper integrals, it is natural to seek
an improper integral corresponding to a power series. If we write the power series as
\begin{eqnarray}\label{gamma1}
\sum_{n=0}^{\infty}a(n)x^{n},
\end{eqnarray}
then a natural analogue is the improper integral
\begin{eqnarray}\label{gamma2}
\int_{0}^{\infty}a(t)x^{t}dt.
\end{eqnarray}
Except for a minor change in notation, this is the \emph{Laplace transform} of the function $a(t)$.

An important particular Laplace transform is the following:
\begin{eqnarray}\label{gamma5}
\int_{0}^{\infty}t^{k}e^{-st}dt,
\end{eqnarray}
in which $a(t)=t^{k}$, the parameter $k$ must be greater than $-1$, to avoid divergence of the integral at $t=0$.
With $k$ so restricted, the integral converges for $s>0$. When $k$ is $0$ or a positive integer, the integral is
easily evaluated:
\begin{eqnarray*}
\int_{0}^{\infty}e^{-st}dt=\frac{1}{s}, \quad \quad \quad \int_{0}^{\infty}te^{-st}dt=\frac{1}{s^{2}},\cdots \quad \quad (s>0).
\end{eqnarray*}
In general, an integration by parts shows that, for $s>0$,
\begin{eqnarray}\label{gamma6}
\int_{0}^{\infty}t^{k}e^{-st}dt=\frac{k}{s}\int_{0}^{\infty}t^{k-1}e^{-st}dt.
\end{eqnarray}
Accordingly, we can apply induction to conclude that, for $s>0$,
\begin{eqnarray}\label{gamma7}
\int_{0}^{\infty}t^{k}e^{-st}dt=\frac{k}{s}\frac{k-1}{s}\cdots\frac{1}{s}\frac{1}{s}=\frac{k!}{s^{k+1}}.
\end{eqnarray}
For $s=1$, (\ref{gamma7}) can be written:

\begin{eqnarray}\label{gamma8}
k!=\int_{0}^{\infty}t^{k}e^{-t}dt.
\end{eqnarray}
This suggest a method of generalizing the factorial, that is, we could use the equation (\ref{gamma8})
to define $k!$ for $k$ an arbitrary real number greater than $-1$. It is customary to denote this generalized
factorial by $\Gamma(k+1)$, the Gamma function $\Gamma(k)$ is defined by the equation
\begin{eqnarray*}
\Gamma(k)=\int_{0}^{\infty}t^{k-1}e^{-t}dt, \quad \quad \quad k>0,
\end{eqnarray*}
when $k$ is positive integer or $0$,
\begin{eqnarray*}
\Gamma(k+1)=k!= \left\{ \begin{array}{ccl}
1.2.\cdots k&, & \quad k>0\\
1&, & \quad k=0.
\end{array} \right.
\end{eqnarray*}

The integral (\ref{gamma5}) is expressible in terms of the Gamma function
\begin{eqnarray*}
\int_{0}^{\infty}t^{k}e^{-st}dt=\frac{\Gamma(k+1)}{s^{k+1}} \quad \quad (s>0).
\end{eqnarray*}
Equation (\ref{gamma6}) then states that
\begin{eqnarray*}
\frac{\Gamma(k+1)}{s^{k+1}}=\frac{k}{s}\frac{\Gamma(k)}{s^{k}},
\end{eqnarray*}
that is,
\begin{eqnarray*}
\Gamma(k+1)=k\Gamma(k).
\end{eqnarray*}
This is the \emph{functional equation of the Gamma function}. The functional equation
can be used to define $\Gamma(k)$ for negative $k$. Thus, we write
\begin{eqnarray*}
\Gamma(\frac{1}{2})=(-\frac{1}{2})\Gamma(-\frac{1}{2}), \quad \quad \quad \Gamma(-\frac{1}{2})=(-\frac{3}{2})\Gamma(-\frac{3}{2}),\cdots
\end{eqnarray*}
in order to define $\Gamma(-1/2), \Gamma(-3/2),\cdots$ in terms of the known value of $\Gamma(1/2)$. This procedure fails only for $k=-1,-2,\cdots$.
In fact, we can show that
\begin{eqnarray*}
\lim_{k\rightarrow\infty}\Gamma(k)=+\infty
\end{eqnarray*}
and, if the Gamma function is extended to negative nonintegral $k$ as above,
\begin{eqnarray*}
\lim_{k\rightarrow -n}|\Gamma(k)|=+\infty
\end{eqnarray*}
for every negative integer $-n$.

Historically, the Gamma function was introduced as
\begin{eqnarray*}
\Gamma(x)=\frac{1}{k} \prod_{n=1}^{\infty} \frac{(1+\frac{1}{n})^{k}}{(1+\frac{k}{n})}
\end{eqnarray*}
by Euler, in 1729. In 1814, Legendre was defined to the Gamma function $\Gamma(x)$ with the improper integral as below:
\begin{eqnarray*}
\Gamma(x)=\int_{0}^{\infty}t^{k-1}e^{-t}dt, \quad  \quad \quad \quad   k>0.
\end{eqnarray*}

\subsection{Almost Convergence}

In theory of sequence spaces, there are many applications of Hahn-Banach Extension Theorem. One of the
most important of these applications is Banach limit which further leads to an important concept of almost
convergence. That is, the $\lim$ functional defined on $c$ can be extended to the whole of $\ell_{\infty}$.
This extended functional is called the \emph{Banach limit}. Lorentz \cite{Lorentz} gave a new type of
convergence definition based on the Banach limit which is called almost convergence, in 1948.\\

Now, we will give some definitions related to the above concept.\\

A linear functional $L$ on $\ell_{\infty}$ is said to be a \emph{Banach limit} if, $L(x)\geq 0$ for $x\geq 0$,
$L(e)=1$, where $e=(1,1,1,\ldots)$ and $L(Sx)=L(x)$, where $S$ is the shift operator defined by $(Sx)_{n}=x_{n+1}$.\\

Let $L$ denotes all Banach limits and $(x_{k})$ be a bounded sequence.A sequence $(x_{k})$ is said to be
\emph{almost convergence} to value $\mathcal{L}$ if all its Banach limits coincide, i.e., $L(x)=\mathcal{L}$. A sequence $(x_{k})$
is said to be \emph{almost $A-$summable} to the value $\mathcal{L}$ if its $A-$transform is almost convergent to $\mathcal{L}$.\\

Lorentz \cite{Lorentz} established the following characterization:\\
A sequence $(x_{k})$ is almost convergent to the number $\mathcal{L}$ if and only if $t_{mn(x)}\rightarrow \mathcal{L}$ as
$m\rightarrow\infty$, uniformly in $n$, where
\begin{eqnarray*}
t_{mn}(x)=\frac{1}{m+1}\sum_{i=0}^{m}x_{n+i}.
\end{eqnarray*}
Let $f$ denote the set of all almost convergent sequences, i.e.,

\begin{eqnarray*}
f=\left\{x=(x_{k})\in \ell_{\infty}: \exists \alpha \in\mathbb{C} \ni \lim_{m\rightarrow\infty}t_{mn}(x)=\mathcal{L}  \quad \textrm{uniformly in } n\right\}.
\end{eqnarray*}
Then, the number $\mathcal{L}$ is called the generalized limit of $x$,
 and we write $\mathcal{L}=f-\lim x_{k}$.\\

There are not essential studies about algebraic structure and topological structure
of the almost convergent space $f$. However, Ba\c{s}ar and Kiri\c{s}ci \cite{BasKir}
have been stated and proved that the space $f$ is non-separable closed subspace of $(\ell_{\infty}, \|.\|_{\infty})$
and the Banach space $f$ has no Schauder basis.\\

We will give some known properties of almost convergent:\\

\begin{itemize}
 \item[(i)] Every convergent sequence is almost convergence at the same limit. But,
  converse of this statement is not true.
    \item[(ii)] The space $f$ is non-separable closed subspace of $(\ell_{\infty}, \|.\|_{\infty})$ \cite{BasKir}. But
    the space $c$ is a separable closed subspace of $(\ell_{\infty}, \|.\|_{\infty})$.
  \item[(iii)] $f$ endowed with the norm $\|.\|_{\infty}$ is a $BK-$space.
      \item[(iv)] $f$ is nowhere dense in $\ell_{\infty}$, dense in itself, and closed and therefore perfect.
  \item[(v)] Almost convergence is not representable by a matrix method.
  \item[(vi)] The $\xi-$dual of $f$ is $\ell_{1}$, where $\xi\in \{\alpha, \beta, \gamma\}$.
\end{itemize}

Following Lorentz \cite{Lorentz}, many mathematicians have considered
almost convergence as summability methods defined by sequences of infinite matrices. King \cite{King}
used the idea of almost convergence to study the almost conservative and almost regular matrices.
In \cite{EiLa}, Eizen and Laush considered the class of almost coercive matrices. Schaefer \cite{Sch}
showed that the classes of almost regular and almost coercive matrices are disjoint. Duran \cite{Duran}
studied the classes of almost strongly regular matrices.

\subsection{Fractional Order Difference Operator}

The idea of the difference of fractional order operator firstly was used by Chapman \cite{Chap} and has been studied
by many mathematicians \cite{And2, And3, Bali, BaliDut2, BaliDut3, DemDuy, BaliDut, Kadak, KadBali, Kutt1, Kutt2}.

Let $(x_{n})$ be any sequence of complex numbers. If $s$ is any real constant, then
\begin{eqnarray*}
\Delta^{s}x_{n}=\sum_{m=0}^{\infty}\binom{m-s-1}{m}x_{n+m},
\end{eqnarray*}
where
\begin{eqnarray*}
\binom{u}{m}=\frac{u(u-1)\cdots(u-m+1)}{m!}
\end{eqnarray*}
and we assume that if $m$ is a negative integer, then we have $\binom{u}{m}=0$.\\

Let $(x_{n})$ be any sequence of complex numbers and $r$, $s$ are positive integers or zero. Then,
\begin{eqnarray}\label{delta1}
\Delta^{r+s}x_{n}=\Delta^{r}(\Delta^{s}x_{n}).
\end{eqnarray}
This equation has been extended to other values of $r$, $s$. The sequence $(x_{n})$  has been mostly chosen as bounded.
If the sequence $(x_{n})$ is bounded then, for $r\geq -1$, $s\geq 0$ and $r+s>0$, the equation (\ref{delta1}) holds
and if $x_{n}\rightarrow 0$ as $n\rightarrow\infty$, then the condition $r+s>0$ may be replaced by $r+s\geq0$ \cite{And1}.
Also, Bosanquet (\cite{Bos1},  Lemma 1) proved the same theorem. Basically, later studies of these type problems have been
in the direction of strengthening the restrictions on $(x_{n})$. It is known that the extending the set of
values in which (\ref{delta1}) is hold. Andersen \cite{And3} was giving a number of these type theorems.
A part of the Andersen's results were proved by Kuttner \cite{Kutt1}. In \cite{Kutt2}, Kuttner mentioned that
the given Theorem A is the best possible result for all concerns the order of magnitude of the individual terms $\Delta^{(r)}x_{n}$.
At the same time, if we choose $s>-1$, then we obtain a stronger result concerning the order of magnitude of their Ces\`{a}ro means.
These new obtained results are best possible in a sense analogues to that of Theorem A.\\

The definition of a sequence of methods of summability of all integer orders, Borel's integral
and exponential methods being particular cases of this sequence has been extended to every real index  $r$, $-\infty<r<\infty$
and investigated to properties of this definition in detail by Ogieveckii \cite{Ogi}. In \cite{Ogi2}, Ogieveckii
obtained new results of earlier study \cite{Ogi}.\\

Baliarsingh \cite{Bali} introduced difference sequence spaces of fractional order and examined some topological
properties and also computed the K\"{o}the-Toeptlitz duals of new spaces. In this work of Baliarsingh is to extended the sequence spaces defined by difference
operator to the spaces of fractional order difference operator.\\

Dutta and Baliarsingh \cite{BaliDut} introduced paranormed difference spaces with fractional order. In \cite{BaliDut3}, Baliarsingh and  Dutta inverstigated some basic properties of the operator $\Delta^{r}$ and gave the applications to the numerical analysis as Newton's forward difference formula and Newton's backward
difference formula. In \cite{KadBali}, Euler sequence spaces are generalized by the operator $\Delta^{r}$. Baliarsingh and  Dutta\cite{BaliDut2} define teh new paranormed spaces derived by fractional order and also Kadak \cite{Kadak} introduced the concept of $\Delta_{v}^{r}-$ strongly Ces\`{a}ro convergence and denote the set of all
$\Delta_{v}^{r}-$ strongly Ces\`{a}ro convergent sequences by $\Delta_{v}^{r}(\omega_{p})$. Demiriz and Duyar \cite{DemDuy} define the new sequence spaces with fractional order difference operator.

\subsection{Sequence spaces}

It is well known that, the $\omega$ denotes the family of all real (or complex)-valued sequences.
$\omega$ is a linear space and each linear subspace of $\omega$ (with the included addition
and scalar multiplication) is called a \emph{sequence space} such as the spaces $c$, $c_{0}$ and
$\ell_{\infty}$, where $c$, $c_{0}$ and $\ell_{\infty}$ denote the set of all convergent
sequences in fields $\mathbb{R}$ or $\mathbb{C}$, the set of all null sequences and the set
of all bounded sequences, respectively. It is clear that the sets $c$, $c_{0}$ and $\ell_{\infty}$
are the subspaces of the $\omega$. Thus, $c$, $c_{0}$ and $\ell_{\infty}$ equipped with a vector space structure,
from a sequence space. By $bs$ and $cs$, we define the spaces of all bounded and convergent series, respectively.\\

The vector space of numerical sequences is called a $K-$ space (or coordinate space). In coordinate space, addition
and scalar multiplication are defined pointwise. An $FK-$ space is a $K-$ space which provided $X$ is a complete
linear metric space. An $FK-$ space whose topology is normable is called a $BK-$ space.\\

Let $A=(a_{nk})$ be an infinite matrix of complex numbers $a_{nk}$ and $x=(x_{k})\in \omega$, where $k,n\in\mathbb{N}$.
Then the sequence $Ax$ is called as the $A-$transform of $x$ defined by the usual matrix product.
Hence, we transform the sequence $x$ into the sequence $Ax=\{(Ax)_{n}\}$, where
\begin{eqnarray}\label{equa1}
(Ax)_{n}=\sum_{k}a_{nk}x_{k}
\end{eqnarray}
for each $n\in\mathbb{N}$, provided the series on the right hand side of (\ref{equa1}) converges for each $n\in\mathbb{N}$.\\

Let $X$ and $Y$ be two sequence spaces. If $Ax$ exists and is in $Y$ for every sequence $x=(x_{k})\in X$, then we say that
$A$ defines a matrix mapping from $X$ into $Y$, and we denote it by writing $A :X \rightarrow Y$ if and only if the series on
the right hand side of (\ref{equa1}) converges for each $n\in\mathbb{N}$ and every $x\in X$, and we have $Ax=\{(Ax)_{n}\}_{n\in \mathbb{N}}\in Y$
for all $x\in X$.  A sequence $x$ is said to be $A$-summable to $l$ if $Ax$ converges to $l$ which is called the $A$-limit of $x$.\\

Let $X$ be a sequence space and $A$ be an infinite matrix. The sequence space
\begin{eqnarray}
X_{A}=\{x=(x_{k})\in\omega:Ax\in X\}
\end{eqnarray}
is called the domain of $A$ in $X$ which is a sequence space.\\

In the book of Basar \cite{Basarkitap}, it can be seen that the qualified studies about matrix domain(also, \cite{BF2}-\cite{BaliDut3}, \cite{candan3}, \cite{candan4}, \cite{candan5}, \cite{kkms}, \cite{kirisci1}-\cite{kirisci4}).


\section{Spaces of almost null and almost convergent sequences of fractional order}

In this section, we will carry to the fractional order operator which is a different difference operator,
to the idea of the difference operators. Using the fractional order operator, we will define the new almost
null and almost convergent sequence spaces and give some properties of new spaces.\\

If we take a positive proper fraction $r$, then we can define the generalized fractional difference operator $\Delta^{(r)}$
as below\cite{Bali, BaliDut2, BaliDut3, BaliDut}:

\begin{eqnarray}\label{mtrx1}
\Delta^{(r)}(x_{k})=\sum_{i=0}^{\infty}(-1)^{i}\frac{\Gamma(r+1)}{i!\Gamma(r-i+1)}x_{k-i}.
\end{eqnarray}

In this text, we suppose that the infinite series \ref{mtrx1} is convergent for all $x\in \omega$, without loss of generality.
The difference operator $\Delta^{(r)}$ can be given by a triangle:
\begin{eqnarray}\label{mtrx2}
\Delta_{nk}^{(r)} = \left\{ \begin{array}{ccl}
(-1)^{n-k}\frac{\Gamma(r+1)}{(n-k)!\Gamma(r-n+k+1)}&, & (0\leq k\leq n),\\
0&, & (k > n).
\end{array} \right.
\end{eqnarray}
\begin{eqnarray*}
\Delta^{(r)}= \left[ \begin{array}{cccccc}
1 & 0 & 0 &0&\ldots \\
-r&1 & 0&0& \ldots \\
\frac{r(r-1)}{2!}&-r & 1 &0 &  \ldots \\
\frac{-r(r-1)(r-2)}{3!}&\frac{r(r-1)}{2!} & -r&1 &  \ldots \\
\vdots & \vdots &  \vdots&  \vdots& \ddots
\end{array} \right].
\end{eqnarray*}

Thus, from the above definition, we can write the inverse of the difference matrix $\Delta_{nk}^{(r)}$ as
\begin{eqnarray*}
\Delta_{nk}^{(-r)} = \left\{ \begin{array}{ccl}
(-1)^{n-k}\frac{\Gamma(1-r)}{(n-k)!\Gamma(-r-n+k+1)}&, & (0\leq k\leq n),\\
0&, & (k > n).
\end{array} \right.
\end{eqnarray*}

\begin{eqnarray*}
\Delta^{(-r)}= \left[ \begin{array}{cccccc}
1 & 0 & 0 &0&\ldots \\
r&1 & 0&0& \ldots \\
\frac{r(r+1)}{2!}&r & 1 &0 &  \ldots \\
\frac{r(r+1)(r+2)}{3!}&\frac{r(r+1)}{2!} & r&1 &  \ldots \\
\vdots & \vdots &  \vdots&  \vdots& \ddots
\end{array} \right].
\end{eqnarray*}

Baliarsingh and Dutta \cite{BaliDut3} proved the following equalities:
\begin{eqnarray}\label{deltaequ}
\Delta^{(r)}\circ \Delta^{(s)}=\Delta^{(r+s)}=\Delta^{(s)}\circ \Delta^{(r)}.
\end{eqnarray}

Let $I_{prop}$ be an identity operator on $\omega$. Then, we have
\begin{eqnarray}\label{deltaequ2}
\Delta^{(r)}\circ \Delta^{(-r)}=\Delta^{(-r)}\circ \Delta^{(r)}=I_{prop}.
\end{eqnarray}

Now, we will define the sequence spaces $fdf$ and $fdf_{0}$ as the set of all sequences
whose $\Delta^{(r)}-$transforms are in the spaces almost convergent sequence and almost null sequence spaces, respectively.

\begin{eqnarray*}
fdf=\left\{x=(x_{k})\in \omega: \exists\mathcal{L}\in\mathbb{C} \ni \lim_{m\rightarrow\infty}t_{mn}\left(\Delta^{(r)}x\right)=\mathcal{L} \quad \textrm{uniformly in } n\right\}
\end{eqnarray*}
and
\begin{eqnarray*}
fdf_{0}=\left\{x=(x_{k})\in \omega: \lim_{m\rightarrow\infty}t_{mn}\left(\Delta^{(r)}x\right)=0 \quad \textrm{uniformly in } n\right\}
\end{eqnarray*}
where
\begin{eqnarray*}
t_{mn}(\Delta^{(r)}x)=\frac{1}{m+1}\sum_{k=0}^{m}\left(\Delta^{(r)}x\right)_{n+k}=
\frac{1}{m+1}\sum_{k=0}^{m}\sum_{j=0}^{n+k}\left[\sum_{i=0}^{k-j}(-1)^{i}\frac{\Gamma(r+1)}{i!\Gamma(r+1-i)}\right]x_{j}
\end{eqnarray*}
uniformly in $n$ and for all $m,n\in \mathbb{N}$. Thus, using the (\ref{mtrx1}), we can write the sequence $y=(y_{k})=\Delta^{(r)}(x_{k})$.

\begin{thm}\label{isomorph}
The sequence spaces $fdf_{0}$ and $fdf$ are linearly isomorphic to the spaces $f_{0}$ and $f$, respectively.
\end{thm}

\begin{proof}
We prove only case $fdf\cong f$. Therefore, we need find a linear bijection between the spaces
$fdf$ and $f$. Now, we define the transformation $T$ such that $T: fdf \rightarrow f$ by $x\mapsto y=Tx=\Delta^{(r)}x$.
It is clear that the transformation $T$ is linear. For, $\theta=(0,0,\ldots)$, if we choose $Tx=\theta$, then,
we obtain $x=\theta$. This shows us that the transformation $T$ is injective.\\

We take the matrix (\ref{mtrx2}) and any $y=(y_{k})\in f$. Then, we can define the sequence $x=(x_{k})$ as
\begin{eqnarray*}
x_{k}=\sum_{j=0}^{\infty}(-1)^{j}\frac{\Gamma(1-r)}{j!(\Gamma(1-r-j))}y_{k-j}-y_{k-j-1}
\end{eqnarray*}
for all $k\in \mathbb{N}$. Consider the equality \label{deltaequ2}. Then, we have,
\begin{eqnarray*}
\lim_{m\rightarrow\infty} \left [t_{mn}\left(\Delta^{(r)}\right)x \right]=
\lim_{m\rightarrow\infty}\frac{1}{m+1}\sum_{i=0}^{m}y_{n+i}=f-\lim y_{n} \quad \textrm{uniformly in } n
\end{eqnarray*}
This means that $x=(x_{k})\in fdf$. It can conclude that the transformation $T$ is surjective and this completes the proof.
\end{proof}

Since the spaces $f_{0}$ and $f$ endowed with the norm $\|.\|_{\infty}$ are $BK-$spaces and
$\Delta^{(r)}$ is a triangle, the spaces $fdf$ and $fdf_{0}$ are $BK-$spaces with the norm $\|.\|_{fdf}$ (\cite{Wil}).
Therefore, there is no need for detailed proof of the following theorem.

\begin{thm}
Spaces of almost null and almost convergent sequences of fractional order $fdf$ and $fdf_{0}$
are $BK-$spaces with the same norm given by
\begin{eqnarray*}
\|x\|_{fdf}=\|\Delta^{(r)}x\|_{f}=\sup_{m,n\in \mathbb{N}}|t_{mn}(\Delta^{(r)}x)|.
\end{eqnarray*}
\end{thm}

In Theorem \ref{inc1} and Theorem \ref{inc2}, we give some inclusion relations between new spaces and
some known spaces.

\begin{thm}\label{inc1}
The sequence spaces $fdf_{0}$ and $fdf$ strictly include the spaces $f_{0}$ and $f$, respectively.
\end{thm}

\begin{proof}
By the definition on the almost null and almost convergent sequence spaces, these inclusions are hold, i.e.,
$f_{0}\subset fdf_{0}$ and $f\subset fdf$.\\

We will show that these inclusions are strict. Therefore, we give the sequence $d=(d_{k})$ defined by
\begin{eqnarray*}
d_{k}=\left\{(-1)^{n-k}\frac{\Gamma(1-r)}{(n-k)!\Gamma(-r-n+k+1)}\left(\frac{1-(-1)^k}{2}\right)\right\}
\end{eqnarray*}
for all $n,k\in \mathbb{N}$. The sequence is not in $f$, but $\Delta^{r}d$ is almost convergent to $1/2$,
as we desired.
\end{proof}

\begin{thm}\label{inc2}
The inclusions $fdf \subset \ell_{\infty}$ and $c\subset fdf$ strictly hold.
\end{thm}

\begin{proof}
It is well known that $c \subset f \subset \ell_{\infty}$. We easily obtain that the inclusions
$c\subset fdf$ and $fdf \subset \ell_{\infty}$ are hold because of Theorem \ref{inc1} and the fact
that $c \subset f \subset \ell_{\infty}$. Further, we must prove that these inclusions are strict.\\

First, we investigate whether the inclusion $fdf \subset \ell_{\infty}$ is strict. Let $r$ be a positive proper fraction.
We choose the sequence
$x=\Delta^{(-r)}y$ with the sequence $y=(y_{k})$ such that the terms of this sequence be it so the blocks
0's are increasing by factors of 100 and the blocks of 1's are increasing by factors of 10(cf. Miller and Orhan \cite{MilOrh}).
This sequence in the set $\ell_{\infty} / f$. Then, the sequence $x$ is not in the space almost convergent sequence of fractional
order but in the space $\ell_{\infty}$, as desired.\\

Second, we must show that the inclusion $c\subset fdf$ is strict. We know that the inclusion $f \subset fdf$ strictly holds from Theorem \ref{inc1}
and the inclusion $c\subset f$ also strictly holds. Therefore, from these results, we can write the inclusion $c\subset fdf$ strictly holds.

\end{proof}

\begin{rem}

It is well known that the almost convergent sequence space $f$ is a non-separable closed subspace of $(\ell_{\infty}, \|.\|_{\infty})$
and as a direct consequence of this case, we can write the Banach space $f$ has no Schauder basis (\cite{BasKir}). For $A=(a_{nk})$ is a triangle and
$X$ is a normed sequence space, the domain $X_{A}$ in a sequence spaces $X$ has a basis if and only if $X$
has a basis.
\end{rem}

Then, we have:

\begin{cor}
The spaces of almost null and almost convergent sequences of fractional order have no Schauder basis.
\end{cor}

If $X,Y\subset \omega$ and $z$ any sequence, we can write $z^{-1}*X=\{x=(x_{k})\in \omega: xz\in X\}$ and $M(X,Y)=\bigcap_{x\in X}x^{-1}*Y$.
If we choose $Y=cs, bs$, then we obtain the $\beta-, \gamma- $duals of $X$, respectively as
\begin{eqnarray*}
X^{\beta}&=&M(X,cs)=\{z=(z_{k})\in \omega:  zx=(z_{k}x_{k})\in cs ~\textrm{for all }~   x\in X\}\\
X^{\gamma}&=&M(X,bs)=\{z=(z_{k})\in \omega:  zx=(z_{k}x_{k})\in bs ~\textrm{for all }~   x\in X\}.
\end{eqnarray*}\\

Now, we will give some lemmas, which are provides convenience in the compute of the dual
spaces and characterize of matrix transformations.

\begin{lem}
Matrix transformations between $BK-$spaces are continuous.
\end{lem}

\begin{lem}\cite[Lemma 5.3]{AB2}\label{mtrxtool0}
Let $X, Y$ be any two sequence spaces. $A\in (X: Y_{T})$ if and only if $TA\in (X:Y)$, where
$A$ an infinite matrix and $T$ a triangle matrix.
\end{lem}

\begin{lem}\cite[Theorem 3.1]{AB}\label{mtrxtool}
We define the $B^{T}=(b_{nk})$ depending on a sequence  $a=(a_{k})\in\omega$ and give the inverse of the triangle matrix $T=(t_{nk})$ by
\begin{eqnarray*}
b_{nk}=\sum_{j=k}^na_{j}v_{jk}
\end{eqnarray*}
for all $k,n\in\mathbb{N}$. Then,
\begin{eqnarray*}
X_{T}^{\beta}=\{a=(a_{k})\in\omega: B^{T}\in(X:c)\}
\end{eqnarray*}
and
\begin{eqnarray*}
X_{T}^{\gamma}=\{a=(a_{k})\in\omega: B^{T}\in(X:\ell_{\infty})\}.
\end{eqnarray*}
\end{lem}

Now, we list the following useful conditions from \cite{bas}, \cite{bslk}, \cite{Duran} and \cite{JAS}.\\
\begin{eqnarray}\label{eq20}
&&\sup_{n\in\mathbb{N}}\sum_{k} |a_{nk}|<\infty,\\\label{eq21}
&&\lim_{n \to \infty}a_{nk}=\alpha_{k}  \quad\quad \textrm{ for each fixed $k\in \mathbb{N}$ }, \\\label{eq22}
&&\lim_{n \to \infty}\sum_{k} a_{nk}=\alpha,\\\label{eq23}
&&\lim_{n \to \infty}\sum_{k} |\Delta(a_{nk}-\alpha_{k})|=0,\\\label{eq24}
&&f-\lim a_{nk}=\alpha_{k}    \quad \textrm{ exists for each fixed $k\in \mathbb{N}$ },\\\label{eq25}
&&f-\lim_{n \to \infty}\sum_{k} a_{nk}=\alpha,\\\label{eq26}
&&\lim_{m \to \infty}\sum_{k} |\Delta(a(n,k,m)-\alpha_{k})|=0 \quad \textrm{uniformly in $n$},\\ \label{eq29}
&&\sup_{n \in\mathbb{N}}\sum_{k}|a(n,k)|<\infty, \\ \label{eq30}
&&\sum_{k}a_{nk}=\alpha_{k} \quad\quad \textrm{ for each fixed $k\in \mathbb{N}$ }, \\ \label{eq31}
&&\sum_{n}\sum_{k}a_{nk}=\alpha, \\ \label{eq32}
&&\lim_{n \to \infty}\sum_{k}\left|\Delta[a(n,k)-\alpha_{k}]\right|=0.
\end{eqnarray}
where $a(n,k,m)=\frac{1}{m+1}\sum_{j=0}^{m}a_{n+j,k}$ and $\Delta a_{nk}=(a_{nk}-a_{n,k+1})$.\\

The following lemma will be used for the Theorem \ref{dualthm} and some examples in Section 4.

\begin{lem}\label{lemMTRX}
For the characterization of the class $(X:Y)$ with
$X=\{f\}$ and $Y=\{\ell_{\infty}, c, cs, bs, f\}$, we can give the necessary and sufficient
conditions from Table 1, where
\begin{itemize}
  \item[\textbf{1.}] (\ref{eq20})
\item[\textbf{2.}] (\ref{eq20}), (\ref{eq21}), (\ref{eq22}), (\ref{eq23})
\item[\textbf{3.}] (\ref{eq29}), (\ref{eq30}), (\ref{eq31}), (\ref{eq32})
\item[\textbf{4.}] (\ref{eq29})
\item[\textbf{5.}] (\ref{eq20}), (\ref{eq24}), (\ref{eq25}), (\ref{eq26})
\end{itemize}
\end{lem}

\begin{center}
\begin{tabular}{|c | c c c c c|}
\hline
To $\rightarrow$ & $\ell_{\infty}$ & $c$ & $cs$ & $bs$ & $f$\\ \hline
From $\downarrow$ &  &  & & &\\ \hline
$f$ & \textbf{1.} & \textbf{2.} & \textbf{3.} & \textbf{4.} & \textbf{5.} \\
\hline
\end{tabular}

\vspace{0.1cm}Table 1\\

\end{center}

In Table 1, we can see that the characterization of the classes $(X:Y)$
with $X=\{f\}$ and $Y=\{\ell_{\infty}, c, cs, bs, f\}$.\\

Let $u_{k}=a_{k}\sum_{i=0}^{k}(-1)^{i}\frac{\Gamma(1-r)}{i!\Gamma(1-r-i)}$. For using in the proof of Theorem \ref{dualthm},
we define the matrix $V=(v_{nk})$ as below:
\begin{eqnarray}\label{mtrxdual}
v_{nk}= \left\{ \begin{array}{ccl}
u_{k}-u_{k+1}&, & (k\leq n),\\
u_{n}&, & (k=n),\\
0&, & (k > n).
\end{array} \right.
\end{eqnarray}

\begin{thm}\label{dualthm}
The $\beta-$ and $\gamma-$ duals of the space $fdf$ defined by
\begin{eqnarray*}
(fdf)^{\beta}=\{u=(u_{k})\in \omega: V\in (f:c)\}
\end{eqnarray*}
and
\begin{eqnarray*}
(fdf)^{\gamma}=\{u=(u_{k})\in \omega: V\in (f:\ell_{\infty})\}
\end{eqnarray*}

\end{thm}

\begin{proof}
Let $u=(u_{k})\in \omega$. We begin the equality
\begin{eqnarray}\label{dual1}
\sum_{k=0}^{n}u_{k}x_{k}=\sum_{k=0}^{n}\left[a_{k}\sum_{i=0}^{\infty}(-1)^{i}\frac{\Gamma(1-r)}{i!\Gamma(1-r-i)}\right](y_{k-i}-y_{k-i-1})
\end{eqnarray}

\begin{eqnarray*}
=\sum_{k=0}^{n-1}\left[a_{k}\sum_{i=0}^{k}(-1)^{i}\frac{\Gamma(1-r)}{i!\Gamma(1-r-i)}-
a_{k+1}\sum_{i=0}^{k+1}(-1)^{i}\frac{\Gamma(1-r)}{i!\Gamma(1-r-i)}\right]y_{k}
\end{eqnarray*}

\begin{eqnarray*}
+\left[a_{n}\sum_{i=0}^{n}(-1)^{i}\frac{\Gamma(1-r)}{i!\Gamma(1-r-i)}\right]y_{n}=(Vy)_{n},
\end{eqnarray*}
where $V=(v_{nk})$ is defined by (\ref{mtrxdual}). Using (\ref{dual1}), we can see that
$ux=(u_{k}x_{k})\in cs$ or $bs$ whenever $x=(x_{k})\in fdf$ if and only if $Vy\in c$ or $\ell_{\infty}$
whenever $y=(y_{k})\in f$. Then, from Lemma \ref{mtrxtool0} and Lemma \ref{mtrxtool}, we obtain  the result that
$u=(u_{k})\in \left(fdf\right)^{\beta}$ or $u=(u_{k})\in \left(fdf\right)^{\gamma}$ if and only if
$V\in (f:c)$ or $V\in (f:\ell_{\infty})$, which is what we wished to prove.
\end{proof}

From this theorem, we have the following corollary:

\begin{cor}\label{dualcor}
The $\beta-$ and $\gamma-$ duals of the space almost convergent sequence of fractional order
\begin{eqnarray*}
fdf^{\beta}= \left\{u=(u_{k})\in \omega: (u_{k}-u_{k+1})\in \ell_{1} \quad \textrm{and} \quad (u_{n})\in c \right\}
\end{eqnarray*}
and
\begin{eqnarray*}
fdf^{\gamma}= \left\{u=(u_{k})\in \omega: (u_{k}-u_{k+1})\in \ell_{1} \quad \textrm{and} \quad (u_{n})\in \ell_{\infty} \right\}.
\end{eqnarray*}
\end{cor}

Therefore, we have defined the new spaces derived by fractional order and examined some structural and topological properties, in this section.
Especially, the $\beta-$ and $\gamma-$ duals of new spaces will help us in the characterization of the matrix transformations.

\section{Some matrix mappings related to the new space}
Let $X$ and $Y$ be arbitrary subsets of $\omega$. We shall show that,
the characterizations of the classes $(X, Y_{T})$ and $(X_{T},Y)$ can be reduced to that of
$(X, Y)$, where $T$ is a triangle.\\

It is well known that if $fdf \cong f$, then the equivalence
\begin{eqnarray*}
x\in fdf \Leftrightarrow y\in f
\end{eqnarray*}
holds. Then, the following theorems will be proved and given
some corollaries which can be obtained to that of Theorems
\ref{mtrxtr1} and \ref{mtrxtr2}. Then, using Lemmas \ref{mtrxtool0} and \ref{mtrxtool}, we have:

\begin{thm}\label{mtrxtr1}
Consider the infinite matrices $A=(a_{nk})$ and $D=(d_{nk})$. These matrices
get associated with each other the following relations:
\begin{eqnarray}\label{eq1}
d_{nk}= \left\{ \begin{array}{ccl}
t_{nk}-t_{n,k+1}&, & (k\leq m),\\
t_{nm}&, & (k=m),\\
0&, & (k > m).
\end{array} \right.
\end{eqnarray}
for all $k,m, n\in \mathbb{N}$, where
\begin{eqnarray*}
t_{nk}=a_{nk}\sum_{j=0}^{k}(-1)^{j}\frac{\Gamma(1-r)}{j!\Gamma(1-i-r)}.
\end{eqnarray*}
Then $A \in (fdf:Y)$ if and only if $\{a_{nk}\}_{k\in\mathbb{N}} \in fdf^{\beta}$
for all $n\in \mathbb{N}$ and $D\in (f:Y)$, where $Y$ be any sequence space.
\end{thm}

\begin{proof}
We assume that the (\ref{eq1}) holds between the entries of $A=(a_{nk})$ and $D=(d_{nk})$.
Let us remember that from Theorem \ref{isomorph}, the spaces $fdf$ and $f$ are linearly isomorphic. Firstly,
we choose any $y=(y_{k})\in f$ and consider $A \in (fdf:Y)$. Then, we obtain that $D\Delta^{r}$ exists and
$\{a_{nk}\}fdf^{\beta}$ for all $k\in \mathbb{N}$. Therefore, the necessity of (\ref{eq1}) yields and
$\{d_{nk}\}\in f^{\beta}$ for all $k,n\in \mathbb{N}$. Hence, $Dy$ exists for each $y\in f$. Thus, if we take
$m\rightarrow \infty$ in the equality
\begin{eqnarray*}
\sum_{k=0}^{m}a_{nk}x_{k}=\sum_{k=0}^{m-1}(t_{nk}-t_{n,k+1})y_{k}+t_{nm}y_{m}
\end{eqnarray*}
for all $m,n\in \mathbb{N}$, then, we understand that $Dy=Ax$. So, we obtain that $D\in (f:Y)$.\\

Now, we consider that $\{a_{nk}\}_{k\in\mathbb{N}} \in fdf^{\beta}$
for all $n\in \mathbb{N}$ and $D\in (f:Y)$. We take any $x=(x_{k})\in fdf$. Then, we can see that
$Ax$ exists. Therefore, for $m\rightarrow \infty$, from the equality

\begin{eqnarray*}
\sum_{k=0}^{m}d_{nk}y_{k}=\sum_{k=0}^{m}\left[\sum_{j=k}^{m}\left(\sum_{i=0}^{j-k}(-1)^{j}\frac{\Gamma(1+r)}{j!\Gamma(1+r-j)}\right)c_{nj}\right]x_{k}
\end{eqnarray*}
for all $n\in \mathbb{N}$, we obtain that $Ax=Dy$. Therefore, this shows that $A \in (fdf:Y)$.
\end{proof}

\begin{thm}\label{mtrxtr2}
Consider that the infinite matrices $A=(a_{nk})$
and $E=(e_{nk})$ with
\begin{eqnarray}\label{eq2}
e_{nk}:=\sum_{j=0}^{n}\left[\sum_{i=0}^{n-j}(-1)^{i}\frac{\Gamma(1+r)}{i!\Gamma(r+1-i)}\right]a_{jk}.
\end{eqnarray}
Then, $A=(a_{nk})\in (X:fdf)$ if and only if $E\in (X :f)$.
\end{thm}

\begin{proof}
We take any $z=(z_{k})\in X$. Using the (\ref{eq2}), we have
\begin{eqnarray}\label{trans1}
\sum_{k=0}^{m}e_{nk}z_{k}=\sum_{j=0}^{n}\left[\sum_{k=0}^{m}\left(\sum_{i=0}^{n-j}(-1)^{i}\frac{\Gamma(1+r)}{i!\Gamma(1+r-i)}\right)a_{jk}\right]z_{k}
\end{eqnarray}
for all $m,n\in \mathbb{N}$. Then, for $m\rightarrow \infty$, equation (\ref{trans1}) gives us that $(Ez)_{n}=\{\Delta^{r}(Az)\}_{n}$.
Therefore, one can immediately observe from this that $Az\in fdf$ whenever $z\in X$ if and only if $Ez\in f$ whenever $z\in X$. Thus, the proof is completed.
\end{proof}

\section{Examples}

In this section, we give some examples related to classes of $A\in (X: fdf)$ and $A\in (fdf:Y)$.\\

If we choose any sequence spaces $X$ and $Y$ in Theorem \ref{mtrxtr1} and \ref{mtrxtr2} in previous section,
then, we can find several consequences in every choice. For example, if we take
the space $\ell_{\infty}$ and the spaces which are isomorphic to $\ell_{\infty}$
instead of $Y$ in Theorem \ref{mtrxtr1}, we obtain the following corollaries and examples:

\begin{cor}\label{corinfinity}
$A\in(fdf:\ell_{\infty})$ if and only if $\{a_{nk}\}\in \{fdf\}^{\beta}$ and
\begin{eqnarray}\label{equamat1}
\sup_{n\in \mathbb{N}}\sum_{k}|d_{nk}|<\infty
\end{eqnarray}
\end{cor}

\begin{cor}\label{corf}
$A\in(fdf:f)$ if and only if $\{a_{nk}\}_{n\in \mathbb{N}}\in \{fdf\}^{\beta}$, (\ref{equamat1}) holds and there are $\alpha_{k},\alpha \in \mathbb{C}$ such that
\begin{eqnarray}\label{eqf10}
&&f-\lim_{n \rightarrow \infty} d_{nk} = \alpha_{k}~ \textrm{ for each }~k\in\mathbb{N}.\\ \label{eqf11}
&&f-\lim_{n \rightarrow \infty} \sum_{k} d_{nk} = \alpha.\\ \label{eqf12}
&&\lim_{n \rightarrow \infty} \sum_{k} |\Delta[d(n,k,m)-\alpha_{k}]| = 0   \textrm{ uniformly in $m$ }.
\end{eqnarray}
where $d(n,k,m)=\frac{1}{n+1}\sum_{j=0}^{n}d_{m+j,k}$.
\end{cor}

\begin{exmp}\label{exmpE}
The Euler sequence space $e_{\infty}^{r}$ is defined by $e_{\infty}^{r}=\{x\in \omega: \sup_{n\in\mathbb{N}}|\sum_{k=0}^{n}\binom{n}{k}(1-r)^{n-k}r^{k}x_{k}|<\infty\}$ (\cite{BF2} and \cite{BFM}).
We consider the infinite matrix $A=(a_{nk})$ and define the matrix $C=(_{nk})$ by
\begin{eqnarray*}
c_{nk}=\sum_{j=0}^{n}\binom{n}{j}(1-r)^{n-j}r^{j}a_{jk}  \quad \quad (k,n\in \mathbb{N}).
\end{eqnarray*}
If we want to get necessary and sufficient conditions for the class $(fdf: e_{\infty}^{r})$ in Theorem \ref{mtrxtr1},
then, we replace the entries of the matrix $A$ by those of the matrix $C$.
\end{exmp}

\begin{exmp}\label{exmpR}
Let $T_{n}=\sum_{k=0}^{n}t_{k}$ and $A=(a_{nk})$ be an infinite matrix. We define the matrix $G=(g_{nk})$ by
\begin{eqnarray*}
g_{nk}=\frac{1}{T_{n}}\sum_{j=0}^{n}t_{j}a_{jk}  \quad \quad (k,n\in \mathbb{N}).
\end{eqnarray*}
Then, the necessary and sufficient conditions in order for $A$ belongs to the class $(fdf:r_{\infty}^{t})$
are obtained from in Theorem \ref{mtrxtr1} by replacing the entries of the matrix $A$ by those of the matrix $G$;
 where $r_{\infty}^{t}$ is the space of all sequences whose $R^{t}-$transforms is in the space $\ell_{\infty}$ \cite{malk}.
\end{exmp}

\begin{exmp}
In the space $r_{\infty}^{t}$, if we take $t=e$, then, this space become to the Cesaro sequence space of non-absolute type $X_{\infty}$ \cite{NgLee}.
As a special case, Example \ref{exmpR} includes the characterization of class $(fdf:r_{\infty}^{t})$.
\end{exmp}

\begin{exmp}
The Taylor sequence space $t_{\infty}^{r}$ is defined by $t_{\infty}^{r}=\{x\in \omega: \sup_{n\in\mathbb{N}}|\sum_{k=n}^{\infty}\binom{k}{n}(1-r)^{n+1}r^{k-n}x_{k}|<\infty\}$ (\cite{kirisci4}).
We consider the infinite matrix $A=(a_{nk})$ and define the matrix $P=(p_{nk})$ by
\begin{eqnarray*}
p_{nk}=\sum_{k=n}^{\infty}\binom{k}{n}(1-r)^{n+1}r^{k-n}a_{jk}  \quad \quad (k,n\in \mathbb{N}).
\end{eqnarray*}
If we want to get necessary and sufficient conditions for the class $(fdf: t_{\infty}^{r})$ in Theorem \ref{mtrxtr1},
then, we replace the entries of the matrix $A$ by those of the matrix $P$.
\end{exmp}

\begin{exmp}
Let $A=(a_{nk})$ be an infinite matrix and the matrix $C=(c_{nk})$ defined by Example \ref{exmpE}. Then, the necessary and sufficient conditions in order for $A$ belongs to the class $(fdf:f(E))$ is obtained from in Corollary \ref{corf} by replacing the entries of the matrix $A$ by those of the matrix $C$; where
\begin{eqnarray*}
f(E)=\{x=(x_{k})\in \omega: \exists l \in \mathbb{C} \ni \lim_{m\rightarrow \infty}\sum_{j=0}^{m}\sum_{k=0}^{n+j}\frac{\binom {n+j}{k}(1-r)^{n+j-k}r^{k}x_{k}}{m+1}=l  \textrm{ uniformly in $n$ }\}
\end{eqnarray*}
defined by Kiri\c{s}ci \cite{kirisci3}.
\end{exmp}

\begin{exmp}
Let $A=(a_{nk})$ be an infinite matrix and the matrix $H=(h_{nk})$ defined by $h_{nk}=sa_{n-1,k}+ra_{nk}$. Then, the necessary and sufficient conditions in order for $A$ belongs to the class $(fdf:\widehat{f})$ is obtained from in Corollary \ref{corf} by replacing the entries of the matrix $A$ by those of the matrix $H$; where
\begin{eqnarray*}
\widehat{f}=\{x=(x_{k})\in \omega: \exists \alpha \in \mathbb{C} \ni \lim_{m\rightarrow \infty}\sum_{j=0}^{m}\frac{sx_{k-1+j}+rx_{k+j}}{m+1}=\alpha  \textrm{ uniformly in $k$ }\}
\end{eqnarray*}
defined by Ba\c{s}ar and Kiri\c{s}ci \cite{BasKir}.
\end{exmp}

\begin{exmp}
Let $A=(a_{nk})$ be an infinite matrix and the matrix $M=(m_{nk})$ defined by $m_{nk}=\sum_{j=k}^{\infty}\frac{a_{nj}}{j+1}$. Then, the necessary and sufficient conditions in order for $A$ belongs to the class $(fdf:\widetilde{f})$ is obtained from in Corollary \ref{corf} by replacing the entries of the matrix $A$ by those of the matrix $M$; where
\begin{eqnarray*}
\widetilde{f}=\{x=(x_{k})\in \omega: \exists \alpha \in \mathbb{C} \ni \lim_{n\rightarrow \infty}\sum_{k=0}^{n}\frac{1}{n+1}\sum_{j=0}^{k}\frac{x_{j+p}}{k+1}=\alpha  \textrm{ uniformly in $p$ }\}
\end{eqnarray*}
defined by Kayaduman and \c{S}eng\"{o}nul \cite{kkms}.
\end{exmp}

If we take the spaces $c$, $cs$ and $bs$ instead of $X$ in Theorem \ref{mtrxtr2}, or $Y$ in Theorem \ref{mtrxtr1}
we can write the following examples. Firstly, we give some conditions and following lemmas:

\begin{eqnarray}\label{eq27}
&&\lim_{n \to \infty}a_{nk}=0  \quad\quad \textrm{ for each fixed $n\in \mathbb{N}$ }, \\\label{eq28}
&&\lim_{n \to \infty}\sum_{k}|\Delta^{2}a_{nk}|=\alpha, \\ \label{eq26x}
&&\lim_{m \to \infty}\sum_{k} |a(n,k,m)-\alpha_{k}|=0 \quad \textrm{uniformly in $n$}, \\ \label{eq33}
&&\sup_{n\in \mathbb{N}}\sum_{k}|\Delta a_{nk}|<\infty
\end{eqnarray}

\begin{lem}
Consider that the $X\in\{\ell_{\infty}, c, bs, cs\}$ and $Y\in \{f\}$.
The necessary and sufficient conditions for $A\in (X:Y)$ can be read the following, from Table 2:
\begin{itemize}
  \item[\textbf{6.}] (\ref{eq20}), (\ref{eq24}), (\ref{eq26x})
\item[\textbf{7.}] (\ref{eq20}), (\ref{eq24}), (\ref{eq25})
\item[\textbf{8.}] (\ref{eq33}), (\ref{eq27}), (\ref{eq24}), (\ref{eq26})
\item[\textbf{9.}] (\ref{eq33}), (\ref{eq24})
\end{itemize}
\end{lem}

\begin{center}
\begin{tabular}{|c | c c  c c|}
\hline
From $\rightarrow$ & $\ell_{\infty}$ & $c$ &  $bs$ & $cs$ \\ \hline
To $\downarrow$ &    & & &\\ \hline
$f$ & \textbf{6.} & \textbf{7.} & \textbf{8.} & \textbf{9.}\\
\hline
\end{tabular}

\vspace{0.1cm}Table 2\\
\end{center}

In Table 2, for the $X\in\{\ell_{\infty}, c, bs, cs\}$ and $Y\in \{f\}$, we give the characterization of the classes $(X : Y )$.\\

\begin{exmp}
We choose $X\in \{fdf\}$ and $Y\in \{\ell_{\infty}, c, cs, bs, f\}$.
The necessary and sufficient conditions for $A\in (X:Y)$ can be taken from the Table 3:
\end{exmp}
\begin{itemize}
  \item[\textbf{1a.}] (\ref{eq20}) holds with $d_{nk}$ instead of $a_{nk}$.
\item[\textbf{2a.}] (\ref{eq20}), (\ref{eq21}), (\ref{eq22}), (\ref{eq23}) hold with $d_{nk}$ instead of $a_{nk}$.
\item[\textbf{3a.}] (\ref{eq29}), (\ref{eq30}), (\ref{eq31}), (\ref{eq32}) hold with $d_{nk}$ instead of $a_{nk}$.
\item[\textbf{4a.}] (\ref{eq29}) holds with $d_{nk}$ instead of $a_{nk}$.
\item[\textbf{5a.}] (\ref{eq20}), (\ref{eq24}), (\ref{eq25}), (\ref{eq26}) hold with $d_{nk}$ instead of $a_{nk}$.
\end{itemize}

\begin{center}
\begin{tabular}{|c | c c c c c|}
\hline
To $\rightarrow$ & $\ell_{\infty}$ & $c$ & $cs$ & $bs$ & $f$\\ \hline
From $\downarrow$ &  &  & & &\\ \hline
$fdf$ & \textbf{1a.} & \textbf{2a.} & \textbf{3a.} & \textbf{4a.} & \textbf{5a.} \\
\hline
\end{tabular}

\vspace{0.1cm}Table 3\\
\end{center}

\begin{exmp}
Consider that the $X\in\{\ell_{\infty}, c, bs, cs\}$ and $Y\in \{fdf\}$.
The necessary and sufficient conditions for $A\in (X:Y)$ can be read the following, from Table 4:
\begin{itemize}
\item[\textbf{6a.}] (\ref{eq20}), (\ref{eq24}), (\ref{eq26x}) hold with $e_{nk}$ instead of $a_{nk}$.
\item[\textbf{7a.}] (\ref{eq20}), (\ref{eq24}), (\ref{eq25}) hold with $e_{nk}$ instead of $a_{nk}$.
\item[\textbf{8a.}] (\ref{eq33}), (\ref{eq27}), (\ref{eq24}) (\ref{eq26}) hold with $e_{nk}$ instead of $a_{nk}$.
\item[\textbf{9a.}] (\ref{eq33}), (\ref{eq24}) hold with $e_{nk}$ instead of $a_{nk}$.
\end{itemize}
\end{exmp}

\begin{center}
\begin{tabular}{|c | c c  c c|}
\hline
From $\rightarrow$ & $\ell_{\infty}$ & $c$ &  $bs$ & $cs$ \\ \hline
To $\downarrow$ &    & & &\\ \hline
$fdf$ & \textbf{6a.} & \textbf{7a.} & \textbf{8a.} & \textbf{9a.}\\
\hline
\end{tabular}

\vspace{0.1cm}Table 4\\
\end{center}

\section{Conclusion}

It is well known that the Gamma Function is an extension of the factorial function and is a very useful tool
in applications. In Gamma function, if $k$ is a positive, then the integral
\begin{eqnarray*}
\int_{0}^{\infty}t^{k}e^{-st}dt=\frac{\Gamma(k+1)}{s^{k+1}} \quad \quad (s>0).
\end{eqnarray*}
converges absolutely. This function is also known as a Euler integral of the second kind. The fractional difference operator
defined by a Gamma function and used some important properties of this function. Therefore, this operator is helping us a lot in practice.\\

The diffrences $\Delta^{r}x_{n}$ were initiated by Chapman \cite{Chap}. Many mathematicians are studied to this operator(cf. \cite{And1}-\cite{And3}, \cite{Bos1}, \cite{Kutt1}), \cite{Kutt2}). Firstly, Baliarsingh \cite{Bali} is introduced to the difference sequence spaces for fractional order. Later, this operator used in the theory of sequence spaces, in several studies(\cite{BaliDut}, \cite{BaliDut2}, \cite{BaliDut3}, \cite{DemDuy}, \cite{Kadak}, \cite{KadBali}).\\

The method of almost convergence has been investigated by Lorentz \cite{Lorentz}. We know that the method of almost convergence is a nonmatrix method.
Some general properties about this method are given by many mathematicians. But, structural properties of the method are not yet obtained.
Despite all this, the method of almost convergence still continues to be important. Recently, following the Ba\c{s}ar and Kiri\c{s}ci \cite{BasKir},
very high quality studies is done(cf. \cite{candan}, \cite{candan2}, \cite{DemKaBas}, \cite{KaraOz}, \cite{KaraOz2}, \cite{kkms}, \cite{kirisci1}, \cite{kirisci2}, \cite{kirisci3}, \cite{murs}, \cite{sonmez}).\\

In this paper, we introduce the method of almost convergence with the fractional order operator and give some topological properties.
We obtain several examples related to this new space. The investigation of the fractional calculus studies of the method of almost convergence
will lead us to new results which are not comparable with the results of this paper.

\section*{Conflict of Interests}
The authors declare that there are no conflict of interests regarding the publication of this paper.

\end{document}